\documentclass[a4paper,12pt,oneside,headsepline]{scrartcl}

\usepackage{etoolbox}
\usepackage{ifdraft}
\usepackage{iftex}

\ifluatex
\usepackage[utf8]{luainputenc}
\else
\usepackage[utf8]{inputenc}
\fi
\usepackage[T1]{fontenc}

\usepackage[english]{babel}

\usepackage{lmodern}
\usepackage{fourier}

\usepackage{amsfonts}
\usepackage{mathrsfs} 
\usepackage{dsfont}

\usepackage{amssymb}

\usepackage{amsmath}
\usepackage{mathtools}
\usepackage[thmmarks, amsmath, amsthm]{ntheorem}
\usepackage{nccmath}    

\usepackage[draft=false]{scrlayer-scrpage}
\usepackage[draft=false]{scrlayer-notecolumn}

\usepackage{needspace}

\providebool{nobiblatex}
\boolfalse{nobiblatex}
\ifbool{nobiblatex}{%
}{%
  \usepackage[backend=biber,style=numeric-comp,giveninits,url=false,sortcites=true,maxbibnames=99]{biblatex}
  \addbibresource{bibliography.bib}
}

\usepackage[final,
            pdfborder={0 0 0}, colorlinks=true,
            linkcolor=BrickRed, citecolor=ForestGreen, urlcolor=RoyalBlue]{hyperref}

\usepackage[cmyk,dvipsnames]{xcolor}

\ifdraft{%
\usepackage[textsize=scriptsize,colorinlistoftodos]{todonotes}
\usepackage{lineno}
\usepackage{scrtime}
}{}

\usepackage{booktabs}
\usepackage[inline]{enumitem}
\usepackage{lettrine}
\usepackage{url}

\KOMAoptions{DIV=12}

\ifdraft{%
\linenumbers
\newcommand*\patchAmsMathEnvironmentForLineno[1]{%
  \expandafter\let\csname old#1\expandafter\endcsname\csname #1\endcsname
  \expandafter\let\csname oldend#1\expandafter\endcsname\csname end#1\endcsname
  \renewenvironment{#1}%
     {\linenomath\csname old#1\endcsname}%
     {\csname oldend#1\endcsname\endlinenomath}}%
\newcommand*\patchBothAmsMathEnvironmentsForLineno[1]{%
  \patchAmsMathEnvironmentForLineno{#1}%
  \patchAmsMathEnvironmentForLineno{#1*}}%
\AtBeginDocument{%
\patchBothAmsMathEnvironmentsForLineno{equation}%
\patchBothAmsMathEnvironmentsForLineno{align}%
\patchBothAmsMathEnvironmentsForLineno{flalign}%
\patchBothAmsMathEnvironmentsForLineno{alignat}%
\patchBothAmsMathEnvironmentsForLineno{gather}%
\patchBothAmsMathEnvironmentsForLineno{multline}%
}}

\clearmainofpairofpagestyles
\automark[section]{subsection}

\renewcommand{\subsectionmark}[1]{}

\lohead{{\small\normalfont \headertitle}}
\rohead{{\small \headerauthors}}

\RedeclareSectionCommand[%
  font=\Large\sffamily\bfseries,%
  beforeskip=1\baselineskip,%
  afterskip=0.5\baselineskip,%
  indent=0em,%
  afterindent=false%
  ]{section}

\RedeclareSectionCommands[%
  font=\normalfont\bfseries,%
  beforeskip=3pt,%
  afterskip=-1em%
  ]{subsection,subsubsection}

\RedeclareSectionCommands[%
  font=\normalfont\itshape,%
  beforeskip=1pt,%
  afterskip=-1em,%
  indent=0pt%
  ]{paragraph}

\ifbool{nobiblatex}{}{%
  \AtEveryBibitem{\clearfield{doi}}
  \AtEveryBibitem{\clearfield{isbn}}
  \AtEveryBibitem{\clearfield{issn}}
  \AtEveryBibitem{\clearfield{pages}}
  \AtEveryBibitem{\clearlist{language}}

  \setlength\bibitemsep{3pt}
  
  \renewbibmacro*{in:}{}
  \DeclareFieldFormat
    [article,inbook,incollection,inproceedings,patent,thesis,unpublished]
    {title}{#1}
}

\newenvironment{enumeratearabic}{
\begin{enumerate}[label=(\arabic*),%
  leftmargin=2.5em,itemindent=0pt,%
  labelindent=.5em,labelwidth=1.5em,labelsep=!,%
  noitemsep]
}{
\end{enumerate}
}

\newenvironment{enumerateroman}{
\begin{enumerate}[label=(\roman*),%
  leftmargin=2.5em,itemindent=0pt,%
  labelindent=.5em,labelwidth=1.5em,labelsep=!,%
  nosep]
}{
\end{enumerate}
}

\newenvironment{enumeratearabic*}{
\begin{enumerate*}[label=(\arabic*)] 
}{
\end{enumerate*}
}

\newenvironment{enumerateroman*}{
\begin{enumerate*}[label=(\roman*)] 
}{
\end{enumerate*}
}

\numberwithin{equation}{section}

\theoremnumbering{arabic}
\newtheorem{theoremcounter}{theoremcounter}[section]
\theoremnumbering{Alph}
\newtheorem{maintheoremcounter}{maintheoremcounter}

\theoremstyle{plain}

\newtheorem{corollary}[theoremcounter]{Corollary}
\newtheorem{lemma}[theoremcounter]{Lemma}

\newtheorem{proposition}[theoremcounter]{Proposition}
\newtheorem{theorem}[theoremcounter]{Theorem}

\theoremstyle{plain}

\newtheorem{maintheorem}[maintheoremcounter]{Theorem}

\theoremstyle{definition}

\newtheorem{example}[theoremcounter]{Example}

\theoremstyle{remark}

\newtheorem{remark}[theoremcounter]{Remark}

\theoremstyle{nonumberremark}

\newtheorem{mainremark}{Remark}

\newenvironment{mainremarkenumerate}
{%
\mainremark
\enumeratearabic
}{%
\endenumeratearabic
\endmainremark
}%

\newcommand{\tx}{\text}

\newcommand{\nbd}{\nobreakdash-\hspace{0pt}}

\makeatletter
\newcommand{\writelabel}[1]{#1\def\@currentlabel{#1}}
\makeatother

\newcommand{\minwidthmathbox}[2]{%
  \mathmakebox[{\ifdim#1<\width\width\else#1\fi}]{#2}%
}

\newcommand{\tbf}{\bfseries}

\newcommand{\bbM}{\ensuremath{\mathbb{M}}}

\newcommand{\rmU}{\ensuremath{\mathrm{U}}}

\newcommand{\ov}{\overline}

\newcommand{\mrelspace}[1]{\mathrel{\mspace{#1}}}

\let\rightarroworig\rightarrow
\renewcommand{\rightarrow}
  {\protect\relbar\mrelspace{-9.7mu}\rightarroworig}

\let\leftarroworig\leftarrow
\renewcommand{\leftarrow}
  {\protect\leftarroworig\mrelspace{-9.7mu}\relbar}

\makeatletter
\@ifundefined{xlongrightarrow}{%

}{%

}
\makeatother

\makeatletter
\@ifundefined{xlongleftarrow}{%

}{%

}
\makeatother

\renewcommand{\Im}{\mathrm{Im}}

\newcommand{\isdiv}{\mathop{\mid}}
\newcommand{\nisdiv}{\mathop{\nmid}}
\newcommand{\isexactdiv}{\mathop{\mid\mspace{-3mu}\mid}}

\renewcommand{\pmod}[1]{\;(\mathrm{mod}\, #1)}

\newenvironment{psmatrix}{\left(\begin{smallmatrix}}{\end{smallmatrix}\right)}

\newcommand{\ZZ}{\ensuremath{\mathbb{Z}}}
\newcommand{\QQ}{\ensuremath{\mathbb{Q}}}

\newcommand{\GL}[1]{\ensuremath{\mathrm{GL}_{#1}}}

\newcommand{\SL}[1]{\ensuremath{\mathrm{SL}_{#1}}}

\renewcommand{\det}{\ensuremath{\mathrm{det}}}

\newcommand{\HS}{\ensuremath{\mathbb{H}}}

\newcommand{\ord}{\operatorname{ord}}
\newcommand{\Ga}{\Gamma}
\newcommand{\ga}{\gamma}
\newcommand{\symdiff}{\mathop{\triangle}}

\newcommand{\bigvert}{\big\vert}

\newcommand{\headertitle}{{%
  Imaginary quadratic fields with $\ell$-torsion-free class groups
}}
\newcommand{\headerauthors}{%
  O.~Beckwith,
  M.~Raum,
  O.~K.~Richter%
}
\title{%
  Imaginary quadratic fields with\\
  $\ell$-torsion-free class groups and\\
  specified split primes
}
\author{%
Olivia Beckwith%
\thanks{The author was partially supported by Simons Foundation Grant~\#953473.}
\and%
Martin Raum%
\thanks{The author was partially supported by Vetenskapsr\aa det Grant~2019-03551.}%
\and%
Olav K. Richter%
\thanks{The author was partially supported by Simons Foundation Grant~\#835652.}
}
\ifdraft{\date{\today\ at\ \thistime}}{\date{}}

\begin{document}
\begingroup
\deffootnote[1em]{1.5em}{1em}{\thefootnotemark}
\maketitle
\endgroup


\begin{abstract}
\small
\noindent
{\tbf Abstract:}
Given an odd prime $\ell$ and finite set of odd primes $S_+$, we prove the existence of an imaginary quadratic field whose class number is indivisible by $\ell$ and which splits at every prime in $S_+$. Notably, we do \emph{not} require that~$p \not\equiv -1 \,\pmod{\ell}$ for any of the split primes~$p$ that we impose. Our theorem is in the spirit of a result by Wiles, but we introduce a new method. It relies on a significant improvement of our earlier work on the classification of non-holomorphic Ramanujan-type congruences for Hurwitz class numbers.
\\[.3\baselineskip]
\noindent
\textsf{\textbf{%
  Class numbers of imaginary quadratic fields%
}}%
\noindent
\ {\tiny$\blacksquare$}\ %
\textsf{\textbf{%
  Hurwitz class numbers%
}}%
\noindent
\ {\tiny$\blacksquare$}\ %
\textsf{\textbf{%
 Ra\-ma\-nu\-jan-type congruences
}}
\\[.2\baselineskip]
\noindent
\textsf{\textbf{%
  MSC Primary:
  11R29%
}}%
\ {\tiny$\blacksquare$}\ %
\textsf{\textbf{%
  MSC Secondary:
  11E41,
  11F33,
  11F37,
  11R11%
}}
\end{abstract}




\Needspace*{4em}
\addcontentsline{toc}{section}{Introduction}
\markright{Introduction}
\lettrine[lines=2,nindent=.2em]{\tbf T}{\,he} presence of~$\ell$-torsion in class groups of number fields has been a rich source of research since the days of Gauss, who developed genus theory to the effect of controlling~$2$-torsion in the class groups of imaginary quadratic fields. Even the next case, $3$\nbd{}torsion is not understood entirely, with landmark contributions by Davenport--Heilbronn~\cite{davenport-heilbroon-1971} in the 70ies which remained the state of the art until Pierce~\cite{pierce-2005} and then Bhargava~\cite{bhargava-2007} in the 2000s accessed the problem through new perspectives (see also~\cite{bhargava-shankar-tsimerman-2013,ellenberg-pierce-wood-2017,pierce-turnage-butterbaugh-wood-2020} for later innovations). They established upper bounds for the size of the $3$-torsion subgroup of the class group of a quadratic field of discriminant~$-D$. The expected bound for the size of~$\ell$-torsion subgroups for any prime~$\ell$ is of the form~$\ll |D|^\epsilon$ for arbitrary~$\epsilon > 0$.

In light of the expectation that~$\ell$-torsion is small, it is natural to ask when it vanishes. This question dates back to Hartung \cite{Hartung}, who showed that there are infinitely many imaginary quadratic fields whose class group has trivial $\ell$-torsion. Much later, Bruinier \cite{Bruinier-1999} used the theory of half integral weight modular forms to show that there exist imaginary quadratic fields with trivial $\ell$-torsion in their class group and which also satisfy any finite set of local conditions at primes which are not congruent to $\pm 1, 0 \pmod{\ell}$. Ono and Skinner \cite{OnoSkinner} proved a similar result to~\cite{Bruinier-1999} for central values of $L$-functions of quadratic twists of newforms. In~2015, Wiles~\cite{wiles-2015} provided an existence theorem of a similar flavor to \cite{Bruinier-1999} (using completely different techniques). More specifically, he showed that there exist imaginary quadratic fields without~$\ell$\nbd{}torsion in their class groups that are split at a finite set of odd primes~$S_+$, inert at a finite set of odd primes~$S_-$, and ramified at a finite set of odd primes~$S_0$. Such a result also has applications to, for instance, canonical periods of elliptic curves (see Corollary~3.4 of~\cite{vatsal-1999}).

Note that the sets~$S_+$, $S_-$, and~$S_0$ in Wiles's theorem cannot be chosen arbitrarily. For example, the primes~$p \in S_+$ must satisfy~$p \not\equiv -1 \,\pmod{\ell}$. Such conditions are typical in the field (see Theorem 7 of \cite{Bruinier-1999}, for example) with the remarkable exception of the results by Bhargava--Shankar--Tsimerman~\cite{bhargava-shankar-tsimerman-2013} and Bhargava--Varma~\cite{bhargava-ila-2016}, which are only available in the case of~$\ell = 3$. The purpose of the present work is to present a new method that likewise allows to discard any conditions on the primes in~$S_+$ but still applies to all primes~$\ell$.  We use that method to establish an existence result in the spirit of~\cite{wiles-2015}.

Throughout, we write~$h(-D)$ for the class number of the quadratic field of discriminant~$-D$.

\begin{maintheorem}%
\label{mainthm:indivisibility}
Let~$\ell \ge 2$ be a prime and~$S_+$ a finite set of odd primes. Then there exists an imaginary quadratic field~$L$ satisfying
\begin{enumerateroman}
\item the class number~$h_L$ of~$L$ is prime to~$\ell$;
\item $L$ is split at each prime of\/~$S_+$.
\end{enumerateroman}
\end{maintheorem}

The analogous theorem in~\cite{wiles-2015} relies on a combination of Galois representations and Shimura curves, while our proof of Theorem~\ref{mainthm:indivisibility} relies on a classification of certain ``non-holomorphic Ramanujan-type congruences'' of Hurwitz class numbers, which is based on the analytic theory of mock modular forms and uses Galois theory only implicitly via a classical result of Serre~\cite{serre-1974,serre-1976}.

We next recall the notion of Ramanujan-type congruences of Hurwitz class numbers and some terminology from our earlier work~\cite{beckwith-raum-richter-2020,beckwith-raum-richter-2022}. Ra\-ma\-nu\-jan-type congruences for the Hurwitz class number~$H(D)$ are congruences of the form
\begin{gather*}
  H(an + b) \equiv 0
  \; \pmod{\ell}
\quad
  \tx{for all\ }n \in \ZZ
\tx{,}
\end{gather*}
where $\ell\geq 3$ is a prime and $a > 0$ and~$b$ are integers (for subtleties concerning the case $\ell=3$, see Section~\ref{sec:extension of previous results}). They generalize the original Ramanujan congruences for the partition function.

Such congruences fall into two different families. On one hand, if~$-b$ is not a square modulo~$a$, then we call them \emph{holomorphic}, since the generating series for $H(an+b)$ is a holomorphic modular form. On the other hand, if~$-b$ is a square modulo~$a$, we call them \emph{non-holomorphic}, since that generating series is a mock modular form, i.e., it has a non-holomorphic modular completion.

One can produce examples via the Hurwitz class number formula in~\eqref{eq:hurwitz_class_number}, which in the holomorphic case are of the form
\begin{gather*}%
 H(5^4 n + 5^3) \equiv 0 \;\pmod{3}
\tx{,}\quad
  H(3^3 n + 3^2) \equiv 0 \;\pmod{5}
\tx{,}\quad
  H(5^3 n + 2 \cdot 5^2) \equiv 0 \;\pmod{7}
\tx{.}
\end{gather*}
In the non-holomorphic case they take the form
\begin{gather*}%
   H(3^3 n + 2\cdot3^2) \equiv 0 \;\pmod{3}
\tx{,}\quad  
  H(5^3 n + 5^2) \equiv 0 \;\pmod{5}
\tx{,}\quad
  H(7^3 n + 3\cdot7^2) \equiv 0 \;\pmod{7}
\tx{.}
\end{gather*}

Two of the notable features of these examples are: In all cases~$a$ is a prime power, and it is a power of~$\ell$ if and only if the congruence is non-holomorphic. This pattern is confirmed by experimental data. Moreover, in previous work we showed that~$\ell \isdiv a, b$ for non-holomorphic congruences, which supports the second observation. More optimistically and equally supported by data, one can hope that~\emph{all} Ramanujan-type congruences of Hurwitz class numbers can be explained by the Hurwitz class number formula~\eqref{eq:hurwitz_class_number}. In this paper, we classify all non-holomorphic Ramanujan-type congruences, and confirm that hope in these cases, under a natural maximality condition which we explain next.

We call a Ra\-ma\-nu\-jan-type congruence for~$H(D)$ modulo~$\ell$ on~$a \ZZ~+~b$ \emph{maximal}, if~$H(D)$ has no Ra\-ma\-nu\-jan-type congruence modulo~$\ell$ on any arithmetic progression~$a' \ZZ~+~b'$ that properly contains~$a \ZZ + b$. Furthermore, we call a Ramanujan-type congruence for~$H(D)$ modulo~$\ell$ on~$a \ZZ~+~b$ \emph{regular} if for all primes~$p \isdiv a$ the~$p$\nbd{}valuation of elements of~$a \ZZ +~b$ is fixed, which is equivalent to~$\ord_p(a) > \ord_p(b)$. We remark that \emph{irregular} (i.e., not regular) Ramanujan-type congruences are equivalent to infinitely many regular ones~(see Proposition~\ref{prop:irregular_atkin_type_congruence}). Finally, we define \emph{maximal regular} Ramanujan-type congruences as those that are maximal among the regular ones. Equivalently, this means that for odd~$p \isdiv a$ we have $\ord_p(a) - 1 = \ord_p(b)$ and if $2 \isdiv a$, then~$\ord_2(a) - 1 \ge \ord_2(b) \ge \ord_2(a) - 3$ (see Theorem~C of our previous work~\cite{beckwith-raum-richter-2022}).

\begin{maintheorem}%
\label{mainthm:classification_nonhol_nondeg}
Let~$\ell \ge 3$ be prime. The maximal regular, non-holomorphic Ramanujan-type congruences of Hurwitz class numbers modulo~$\ell$ occur exactly on~$a \ZZ + b = \ell^{m+1} \ZZ + u \ell^m$ for any positive, even integer~$m$ and any~$u \in \ZZ$ with~$(\frac{-u}{\ell}) = +1$. In particular, we have~$\ell^3 \isdiv a$, $a \isexactdiv \ell b$, and~$(\frac{-\ell b \slash a}{\ell}) = +1$.
\end{maintheorem}

\begin{mainremarkenumerate}
\item
Theorem~\ref{mainthm:classification_nonhol_nondeg} significantly improves our results in~\cite{beckwith-raum-richter-2020,beckwith-raum-richter-2022}, which only showed that if $\ell\geq 5$ is prime and if there is a maximal non-holomorphic Ramanujan-type congruence on Hurwitz class numbers modulo~$\ell$ on~$a \ZZ + b$, then~$\ell \isdiv a, b$ and~$a \isdiv \ell b$ (see also Section~\ref{sec:extension of previous results} for the case~$\ell = 3$).

\item
Zagier's Eisenstein series can be viewed as an analog of mock theta series of weight~$\frac{3}{2}$, and it is natural to ask for results like ours in that setting. A classification is currently not available, but Corollary~1.2 and Theorem~1.4 in~\cite{andersen-2014} provide some analog of our results in~\cite{beckwith-raum-richter-2020,beckwith-raum-richter-2022} for the mock theta functions~$f(q)$ and~$\omega(q)$ on~$a \ZZ + b + \frac{1}{24}$ and~$a \ZZ + b - \frac{2}{3}$, respectively. However, in contrast to the situation of Hurwitz class numbers, there are no known examples of non-holomorphic Ramanujan-type congruences for~$f(q)$ or~$\omega(q)$.
\end{mainremarkenumerate}

The connection between Theorem~\ref{mainthm:indivisibility} and~\ref{mainthm:classification_nonhol_nondeg} is made by the following result.

\begin{maintheorem}
\label{mainthm:ramanujan_type_and_discriminant_conditions}
Let~$\ell \ge 3$ be prime.  Assume that there is a maximal regular Ramanujan-type congruence for Hurwitz class numbers modulo~$\ell$ on~$a \ZZ + b$ for cube-free, odd~$a$. Define a disjoint union of finite sets of odd primes~$S = S_+ \cup S_- \cup S_0$ via the formulae
\begin{gather*}
  S_+ := \big\{ p \,:\,  p \isexactdiv a, \big(\tfrac{-b}{p}\big) = +1 \big\}
\tx{,}\quad
  S_- := \big\{  p \,:\, p \isexactdiv a, \big(\tfrac{-b}{p}\big) = -1 \big\}
\tx{,}\quad
  S_0 := \big\{  p \,:\, p^2 \isexactdiv a, \big(\tfrac{-b}{p}\big) = 0 \big\}
\tx{,}
\end{gather*}
where~$p$ always denotes an odd prime. Then for any imaginary quadratic field~$L$ that
\begin{enumerateroman}
\item is split at each prime of\/~$S_+$,
\item is inert at each prime of\/~$S_-$,
\item is ramified at each prime of\/~$S_0$,
\end{enumerateroman}
the class number~$h_L$ is divisible by~$\ell$.

Conversely, if the conclusion in the previous paragraph holds for some disjoint union of finite sets of odd primes~$S = S_+ \cup S_- \cup S_0$, then there is a maximal regular Ramanujan-type congruence for Hurwitz class numbers modulo~$\ell$ on all~$a \ZZ + b$ with
\begin{gather*}
  a = \prod_{p \in S_+ \cup S_-} p\,
      \prod_{p \in S_0} p^2
\tx{,}\qquad
  b \in \ZZ \tx{\ such that\ }
  \big( \tfrac{-b}{p} \big) = \pm 1 \tx{\ if\ } p \in S_\pm
  \tx{\ and\ }
  p \isexactdiv b \tx{\ if\ } p \in S_0
\tx{.}
\end{gather*}
\end{maintheorem}

\begin{mainremark}
The Ramanujan-type congruences in Theorem~\ref{mainthm:ramanujan_type_and_discriminant_conditions} are non-holomorphic if and only if~$S_- = S_0 = \emptyset$. This and the restriction to the non-holomorphic case in Theorem~\ref{mainthm:classification_nonhol_nondeg} explain the absence of~$S_-$ and~$S_0$ in Theorem~\ref{mainthm:indivisibility} compared to Wiles's result~\cite{wiles-2015}. Moreover, the Ramanujan-type congruences on~$a \ZZ + b$ in Theorem~\ref{mainthm:ramanujan_type_and_discriminant_conditions} have square-free~$a$ if~$S_- = S_0 = \emptyset$, which by Theorem~\ref{mainthm:classification_nonhol_nondeg} cannot exist. Thus, the equivalence in Theorem~\ref{mainthm:ramanujan_type_and_discriminant_conditions} relates the impossibility of such congruences to the existence of imaginary quadratic fields with~$\ell$\nbd{}torsion-free class groups as in Theorem~\ref{mainthm:indivisibility}.
\end{mainremark}

The paper is organized as follows.  In Section~\ref{sec:extension of previous results}, we extend our previous results on Ra\-ma\-nu\-jan-type congruences for Hurwitz class numbers modulo~$\ell\geq 5$ to include the case $\ell=3$. In Section~\ref{sec:atkin_type_congruence}, we introduce Atkin-type congruences for Hurwitz class numbers. In Section~\ref{sec:modular_substitution_formula_Radu}, we provide a modular substitution formula \`a la~\cite{radu-2013} that is needed for our proof of Theorem~\ref{mainthm:ramanujan_type_and_discriminant_conditions}. In Section~\ref{sec:Proof_of_Theorem_B}, we prove Theorem~\ref{mainthm:ramanujan_type_and_discriminant_conditions}. Finally, in Section~\ref{sec:Proof_of_Theorems_A_C} we establish Theorems~\ref{mainthm:indivisibility} and~\ref{mainthm:classification_nonhol_nondeg}.


\section{Extension of previous results on Ramanujan-type congruences for Hurwitz class numbers}
\label{sec:extension of previous results}

As in the introduction, $H(D)$ denotes the Hurwitz class number and $h(-D)$ stands for the class number of the quadratic field of discriminant~$-D$. In particular, if $-D<-4$ is a fundamental discriminant, then $H(D)=h(-D)$. Moreover, write~$\omega(-D f^2)$ for the number of units in the quadratic order of discriminant~$-D f^2$. We now state the Hurwitz class number formula, which we use frequently in this work.
\begin{gather}
\label{eq:hurwitz_class_number}
\begin{aligned}
  H(D f^2)
&{}=
  H(D)
  \frac{\omega(-D f^2)}{\omega(-D)}\,
  \sum_{d \isdiv f} d
  \prod_{\substack{p \isdiv d\\p\tx{\ prime}}} \!\!
  \big( 1 - \tfrac{1}{p} \big(\tfrac{-D}{p}\big) \big)
\\
&{}=
  H(D)
  \frac{\omega(-D f^2)}{\omega(-D)}
  \prod_{\substack{p \isdiv f\\p\tx{\ prime}}} \!\!
  \Big(
  1 \,+\,
  \sum_{\substack{d \isdiv f\\ d \tx{\ a\ }p\tx{-power}\\ d \ne 1}} \mspace{-20mu}
  \big( d - \tfrac{d}{p} \big(\tfrac{-D}{p}\big) \big)
  \Big)
\tx{.}
\end{aligned}
\end{gather}

Recall from our earlier work~\cite{beckwith-raum-richter-2020,beckwith-raum-richter-2022} the following Theorems:

\begin{theorem}%
\label{thm:ell_div_ab}
Let~$\ell \ge 5$ be prime. If~$a \ZZ + b$ supports a non-holomorphic Ramanujan-type congruence for Hurwitz class numbers modulo~$\ell$, then~$\ell \isdiv a, b$.
\end{theorem}

\begin{theorem}%
\label{thm:rc_on_square_classes}
Let~$\ell \ge 5$ be prime. If~$a \ZZ + b$ supports a Ramanujan-type congruence for Hurwitz class numbers modulo~$\ell$, then it has a Ramanujan-type congruence on~$a \ZZ + u^2 b$ for all~$u \in \ZZ$, $\gcd(a,u) = 1$.
\end{theorem}

Observe that our papers~\cite{beckwith-raum-richter-2020,beckwith-raum-richter-2022} assumed that~$\ell \ge 5$. The primary purpose of this assumption was to exclude the case of fractional class numbers, in particular the cases~$H(3) = \frac{1}{3}$ and~$H(4) = \frac{1}{2}$.  However, one can extend congruences modulo~$\ell$ to all rational numbers by~$x \equiv y \,\pmod{\ell}$ if and only if~$x-y$ is~$\ell$\nbd{}integral and~$\ell \isdiv (x-y)$.  Then Ramanujan-type congruences of Hurwitz class numbers modulo any prime~$\ell$ can be defined in the usual way:
\begin{gather}
\label{eq:def:ramanujan_type_congruence_ell3}
  \forall n \in \ZZ \,:\,
  H(a n + b) \equiv 0
  \;\pmod{\ell}
\tx{.}
\end{gather}

The purpose of this section is to extend Theorems~\ref{thm:ell_div_ab} and~\ref{thm:rc_on_square_classes} to include the case $\ell=3$.

\begin{theorem}%
\label{thm:ell_div_ab_ell3}
Let~$\ell \ge 3$ be prime. If~$a \ZZ + b$ supports a non-holomorphic Ramanujan-type congruence for Hurwitz class numbers modulo~$\ell$, then~$\ell \isdiv a, b$.
\end{theorem}

\begin{theorem}%
\label{thm:rc_on_square_classes_ell3}
Let~$\ell \ge 3$ be prime. If~$a \ZZ + b$ supports a Ramanujan-type congruence for Hurwitz class numbers modulo~$\ell$, then it has a Ramanujan-type congruence on~$a \ZZ + u^2 b$ for all~$u \in \ZZ$, $\gcd(a,u) = 1$.
\end{theorem}

\begin{remark}
One can also study Ramanujan-type congruences of Hurwitz class numbers modulo~$2$. The Hurwitz class number formula~\eqref{eq:hurwitz_class_number}  together with Gauss's genus theory yields a characterization of such congruences. Let~$a > 0$ and~$b$ be integers. Assume without loss of generality that~$\ord_p(b) = \ord_p(a) - 1$ for odd primes~$p \isdiv a$ and~$\ord_2(b) = \ord_2(a) - 3$ if\/~$2 \isdiv a$. Then we have a Ramanujan-type congruence of the Hurwitz class numbers modulo~$2$ on~$a \ZZ + b$ if and only if neither of the following conditions hold:
\begin{enumerateroman}
\item $p^{1 + 4 \ZZ} \equiv b \,\pmod{a}$ for some odd prime~$p \equiv -1 \,\pmod{4}$ or
\item $4 f^2 \equiv b \,\pmod{a}$ for an integer~$f$.
\end{enumerateroman}
\end{remark}

The proofs of Theorems~\ref{thm:ell_div_ab_ell3} and~\ref{thm:rc_on_square_classes_ell3} are completely analogous to the proofs of Theorems~\ref{thm:ell_div_ab} and~\ref{thm:rc_on_square_classes}.  Recall that our work~\cite{beckwith-raum-richter-2020,beckwith-raum-richter-2022} relied on the following result by Serre, which in~\cite{beckwith-raum-richter-2020,beckwith-raum-richter-2022} we stated for primes $\ell\geq 5$, but which also holds for $\ell=3$. Here $\Ga_1(N)$ denotes the usual congruence subgroup of~$\SL{2}(\ZZ)$.

\begin{theorem}[{Deligne-Serre~\cite{deligne-serre-1974} and Serre~\cite{serre-1974,serre-1976}}]
\label{thm:serre}
Fix a prime~$\ell\geq 3$ and~$k, N \in \ZZ_{\ge 1}$.\linebreak Then there are infinitely many primes~$p \equiv 1 \,\pmod{\ell N}$ such that for every modular form~$f$ of weight~$k$ for~$\Ga_1(N)$ with $\ell$-integral Fourier coefficients, we have
\begin{gather}
\label{eq:thm:serre}
  c(f;\, n p^r)
\;\equiv\;
  (r + 1) c(f;\, n)
  \;\pmod{\ell}
\end{gather}
for all $n \in \ZZ$ co-prime to~$p$ and all~$r \in \ZZ_{\ge 0}$.
\end{theorem}

The extension of Serre's result to quasimodular forms needs a slight modification for the case $\ell~=~3$. See~\cite{zagier-1994,kaneko-zagier-1995} for a background on quasimodular forms.

\begin{corollary}
\label{cor:serre}
Fix a prime~$\ell\geq 3$ and positive integers~$k$ and~$N$. Then there are infinitely many primes~$p \equiv 1 \,\pmod{\ell N}$ such that for every quasi-modular form~$f$ of weight~$k$ for~$\Ga_1(N)$ with $\ell$-integral Fourier coefficients, we have the congruence~\eqref{eq:thm:serre}.
\end{corollary}
\begin{proof}
If $\ell\geq 5$, then Corollary~\ref{cor:serre} reduces to Corollary~1.2 of our paper~\cite{beckwith-raum-richter-2020}, and it suffices to justify the case $\ell=3$.  Recall the weight~$k$ Eisenstein series
\begin{gather*}
  E_k(\tau)
\;:=\;
  1 - \frac{2k}{B_k} \sum_{n = 1}^\infty \sigma_{k-1}(n) \exp(2\pi i\, n\tau)
\tx{,}\quad
  \sigma_{k-1}(n)
\;:=\;
  \sum_{\substack{d \isdiv n\\d>0}} d^{k-1}
\tx{,}
\end{gather*}
where $E_2$ is a quasi-modular form and $E_k$ is a modular form if $k\geq 4$ is even.  Moreover, recall that every quasi-modular form of weight~$k$ for~$\Ga_1(N)$ can be written as
\begin{gather*}
  \sum_{n = 0}^d
  E_2^n f_n
\tx{,}
\end{gather*}
where $f_n$ are modular forms of weights $k-2n$ for~$\Ga_1(N)$.  Note that $E_2\equiv E_4\equiv E_6\equiv 1\pmod{3}$. Set $a(n)=\frac{1+2n+3(-1)^n}{4}$ and $b(n)=\frac{1-(-1)^n}{2}$.  Then
\begin{gather*}
  \sum_{n = 0}^d
  E_2^n f_n
\;\equiv\;
  \sum_{n = 0}^d
  E_4^{a(n)}E_6^{b(n)}f_n 
  \;\pmod{3}
\tx{,}
\end{gather*}
and applying Theorem~\ref{thm:serre} to the modular form of weight~$k + 4$ on the right hand side yields the claim for $\ell=3$.
\end{proof}

Now one can proceed exactly as in~\cite{beckwith-raum-richter-2020,beckwith-raum-richter-2022} to prove Theorems~\ref{thm:ell_div_ab_ell3} and~\ref{thm:rc_on_square_classes_ell3}, and the applications of Corollary~\ref{cor:serre} do not cause any problems if $\ell=3$. Specifically, the proof of these statements is based on calculations that are~$\ell$-integral and  the only consequence of having a Ramanujan-type congruence that these proofs use is that the holomorphic part of~$\rmU_{a,b}\, E_{\frac{3}{2}}$ that appears there is congruent to~$0$ modulo $\ell$, which remains true for $\ell=3$. For brevity, we do not repeat the details.

\section{Atkin-type congruences for Hurwitz class numbers}
\label{sec:atkin_type_congruence}

Theorems~\ref{thm:rc_on_square_classes} and~\ref{thm:rc_on_square_classes_ell3} in Section~\ref{sec:extension of previous results} lead to the notion of multiplicative congruences in the sense of Atkin~\cite{atkin-1968}, which we will call \emph{ Atkin-type congruences}. In order to give a formal definition, we introduce some more notation.  

Given a nonzero integer~$n$ and a prime~$p$, we throughout write~$n = n_p n_p^\#$ for its factorization into a~$p$-power~$n_p$ and an integer~$n_p^\#$ that is co-prime to~$p$. If~$n = 0$, then set~$n_p := 1$ and~$n_p^\# :=~0$.  If~$-n$ is a discriminant, then one has the factorization~$-n = -D f^2 \in \ZZ_{< 0}$, where $-D$ is a fundamental discriminant and $f$~is a positive integer. When writing~$-D f^2$ for a negative integer, we will always assume that~$-D$ is a fundamental discriminant and $f$~is a positive integer.  For integers~$a > 0$ and~$b$ set
\begin{gather}
  b\, (\ZZ \slash a \ZZ)^{\times\, 2}
:=
  \big\{
  n \in \ZZ \,:\,
  \exists u \in \ZZ, \gcd(a,u) = 1 ,\,
  n \equiv b u^2 \,\pmod{a}
  \big\}
\tx{.}
\end{gather}
Note that~$b (\ZZ \slash a \ZZ)^{\times 2}$ is determined by the valuation of~$b$ at primes~$p \isdiv a$ and by the square-class of~$b_p^\#$ modulo~$p$ if~$p$ is odd, or modulo~$\gcd(8,a)$ if~$p = 2$.

We say that Hurwitz class numbers have an \emph{Atkin-type congruence modulo~$\ell$ around $b \pmod{a}$} if
\begin{gather}
\label{eq:def:atkin_type_congruence}
  \forall n \in b (\ZZ \slash a \ZZ)^{\times\, 2} \,:\,
  H(n) \equiv 0 \,\pmod{\ell}
\tx{.}
\end{gather}

As a direct consequence of Theorem~\ref{thm:rc_on_square_classes_ell3} and the inclusion~$a \ZZ + b \in b (\ZZ \slash a \ZZ)^{\times\,2}$, we have the following characterization.

\begin{corollary}%
\label{cor:aktin_type_ramanujan_type}
The Hurwitz class numbers have a Ramanujan-type congruence modulo a prime~$\ell \ge 3$ on the arithmetic progression~$a \ZZ + b$ if and only if they have an Atkin-type congruence modulo~$\ell$ around~$b \pmod{a}$.
\end{corollary}

The notions of maximal, regular, and non-holomorphic Ramanujan-type congruences translate to Atkin-type congruences. More precisely, an Atkin-type congruence around $b \pmod{a}$ is \emph{maximal}, if there is none on any~$b' \pmod{a'}$ with~$b' (\ZZ \slash a' \ZZ)^{\times\,2} \subsetneq b (\ZZ \slash a \ZZ)^{\times\,2}$. It is \emph{regular at~$p \isdiv a$}, if the $p$\nbd{}valuation of elements of~$b (\ZZ \slash a \ZZ)^{\times\,2}$ is fixed, that is, if~$\ord_p(a) > \ord_p(b)$, and \emph{regular} if it is regular for all~$p \isdiv a$. Finally, it is \emph{non-holomorphic} if~$-b$ is a square modulo~$a$, that is, if~$-b (\ZZ \slash a \ZZ)^{\times\,2}$ contains a square. We now give the translation of Theorem~\ref{thm:ell_div_ab_ell3} for Atkin-type congruences.

\begin{theorem}%
[Reformulation of Theorem~\ref{thm:ell_div_ab_ell3}]%
\label{thm:ell_div_ab_atkin_type}
Let~$\ell \ge 3$ be prime.  Assume that Hurwitz class numbers have a non-holomorphic Atkin-type congruence modulo~$\ell$ around~$b \pmod{a}$. Then~$\ell \isdiv a, b$.
\end{theorem}

The proof of Theorem~C in~\cite{beckwith-raum-richter-2022} on Ramanujan-type congruences extends also to Atkin-type congruences, and we now record that extension without repeating the proof.

\begin{proposition}%
\label{prop:maximal_atkin_type_ordp}
Let~$\ell \ge 3$ be prime.  Assume that Hurwitz class numbers have a maximal Atkin-type congruence modulo~$\ell$ around~$b \pmod{a}$. Then for odd primes~$p$,
\begin{gather*}
  \ord_p\big( a \slash \gcd(a,b) \big) \le 1
\quad\tx{and}\quad
  \ord_2\big( a \slash \gcd(a,b) \big) \le 3
\tx{.}
\end{gather*}
\end{proposition}

The next proposition clarifies the role of irregular Atkin-type congruences in the context of this paper.  We further illustrate their role in the example that follows.

\begin{proposition}%
\label{prop:irregular_atkin_type_congruence}
Let~$\ell \ge 3$ be prime. Assume that Hurwitz class numbers have an Atkin-type congruence modulo~$\ell$ around~$b \pmod{a}$ that is irregular at a prime~$p$. Without loss of generality assume that~$a_p \isexactdiv b$, which can be achieved by replacing~$b$ by~$b + a$ if necessary. Then we have regular-at-$p$ Atkin-type congruences of the Hurwitz class numbers modulo~$\ell$ around
\begin{gather*}
  b u p^m\, \pmod{a p^{m+1}}
\quad
 \tx{for all\ } m \in \ZZ_{\ge 0}
 \tx{\ and all\ } u \in \ZZ 
 \tx{\ with\ } p \nisdiv u \tx{\ and\ } u p^m \equiv 1 \,\pmod{a_p^\#}
\tx{.}
\end{gather*}
In particular, non-holomorphic, irregular Atkin-type congruences imply both holomorphic and non-holomorphic ones that are regular.

Conversely, if we have Atkin-type congruences of the Hurwitz class numbers modulo~$\ell$ around all such~$b u p^m \,\pmod{a p^{m+1}}$, then we have an irregular one around~$b \,\pmod{a}$.
\end{proposition}


\begin{proof}
One direction follows directly from the inclusions
\begin{gather*}
  b u p^m \big( \ZZ \slash a p^{m+1} \ZZ \big)^{\times\,2}
\subseteq
  b u p^m \big( \ZZ \slash a \ZZ \big)^{\times\,2}
\subseteq
  b \big( \ZZ \slash a \ZZ \big)^{\times\,2}
\tx{,}
\end{gather*}
where the second inclusion follows from the Chinese Remainder Theorem and the congruence~$b u p^m \equiv 0 \equiv b \,\pmod{a_p}$. The other direction is a consequence of the fact that
\begin{gather*}
  b \big( \ZZ \slash a \ZZ \big)^{\times\,2}
=
  \bigcup_{u, m}
  b u p^m \big( \ZZ \slash a p^{m+1} \ZZ \big)^{\times\,2}
\tx{.}
\end{gather*}
\end{proof}

\begin{example}
We apply Proposition~\ref{prop:irregular_atkin_type_congruence}, Theorem~\ref{mainthm:ramanujan_type_and_discriminant_conditions} (proved in  Section~\ref{sec:Proof_of_Theorem_B}), and the main result of Wiles~\cite{wiles-2015} (not used anywhere else in our paper) to show that there are no Atkin-type congruences for Hurwitz class numbers modulo primes~$\ell \ge 3$ around~$0\, (\ZZ \slash \ell^2 \ZZ)^{\times, \, 2}$.  

Suppose there was an Atkin-type congruence for Hurwitz class numbers modulo~$\ell \ge 3$ around~$0\, (\ZZ \slash \ell^2 \ZZ)^{2\,\times}$, then (after replacing $b=0$ with $b+a=\ell^2$) Proposition~\ref{prop:irregular_atkin_type_congruence} implies for example the following four regular Atkin-type congruences for Hurwitz class numbers modulo~$\ell \ge 3$ around:
\begin{gather*}
  -\ell^2 \pmod{\ell^3}
\tx{,}\quad
  - u \ell^2 \pmod{\ell^3}
\tx{,}\quad
  -\ell^3 \pmod{\ell^4}
\tx{,}\quad
  - u \ell^3 \pmod{\ell^4}
\tx{,}
\end{gather*}
where~$u$ is an integer such that~$(\frac{u}{\ell}) = -1$. The first congruence is non-holomorphic, and follows directly from the Hurwitz class number formula~\eqref{eq:hurwitz_class_number}. The later three congruences are holomorphic, and an application of the Hurwitz class number formula~\eqref{eq:hurwitz_class_number} yields the following maximal Atkin-type congruence for Hurwitz class numbers modulo~$\ell \ge 3$ around
\begin{gather*}
  - u \pmod{\ell}
\tx{,}\quad
  -\ell \pmod{\ell^2}
\tx{,}\quad
  - u \ell \pmod{\ell^2}
\tx{.}
\end{gather*}
Observe here that if~$D f^2$ in~$-u \ell^2 (\ZZ \slash \ell^3 \ZZ)$, in~$-\ell^3 (\ZZ \slash \ell^4 \ZZ)$, and in~$-u \ell^3 (\ZZ \slash \ell^4 \ZZ)$, respectively, then the contribution of~$p = \ell$ in the product over~$p \isdiv f$ in~\eqref{eq:hurwitz_class_number} equals~$1 + (\ell + 1)$, $1 + \ell$, and~$1 + \ell$, respectively. Now consider sets $S_+$, $S_-$, and $S_0$ as in Theorem~\ref{mainthm:ramanujan_type_and_discriminant_conditions}. Then the main result of Wiles~\cite{wiles-2015} guarantees the existence of an imaginary quadratic field~$L$ such that the class number~$h_L$ of~$L$ is prime to~$\ell$ which contradicts Theorem~\ref{mainthm:ramanujan_type_and_discriminant_conditions}.  Thus, there is no Atkin-type congruence modulo~$\ell$ around~$0 \pmod{\ell^2}$.
\end{example}

\section{A modular substitution formula \`{a} la Radu}%
\label{sec:modular_substitution_formula_Radu}

In this section we establish a harmonic Maass form analog of Lemma 5.1 of~\cite{radu-2013}, which is required for our proof of Theorem~\ref{thm:congruences_coupled_square_classes}. We first introduce some necessary standard notation. Throughout, $\tau\in\HS$ (the usual complex upper half plane) and $e(s\tau) := \exp(2 \pi i\, s\tau)$ for $s\in\QQ$. For odd~$D$, set
\begin{gather*}
  \epsilon_D
=
  \begin{cases}
  1\tx{,} & \tx{if $D \equiv 1 \,\pmod{4}$\tx{;}} \\
  i\tx{,} & \tx{if $D \equiv 3 \,\pmod{4}$\tx{.}}
  \end{cases}
\end{gather*}
For $\ga = \begin{psmatrix} a & b \\ c & d \end{psmatrix} \in \GL{2}(\QQ)$ with~$\det(\ga) > 0$, the weight-$k$ slash operator is defined by
\begin{gather}
\label{eq:slash_k}
  \big( f \bigvert_k \ga \big) (\tau)
\;=\;
  (\det\,\ga)^{\frac{k}{2}}\,
  (c \tau + d)^{-k}\,
  f \big(\mfrac{a \tau + b}{c \tau + d}\big)
\tx{.}
\end{gather}

Let $\bbM_k(\Ga_0(N))$ be the space of harmonic Maass forms (satisfying the moderate growth condition at all cusps) of half-integral weight~$k$ for the usual congruence subgroup~$\Ga_0(N) \subseteq \SL{2}(\ZZ)$ with respect to the multiplier $\nu_{\theta}^{2k}$, where $\nu_{\theta}$ is the theta multiplier. For more details on modular forms and harmonic Maass forms, see for example~\cite{lang-1995} and~\cite{bringmann-folsom-ono-rolen-2018}.

Finally, for $a \in \ZZ_{\ge 1}$ and $b \in \ZZ$ recall the operators~$\rmU_{a,b}$ from our earlier work~\cite{beckwith-raum-richter-2020,beckwith-raum-richter-2022}, which act on Fourier series expansions of non-holomorphic modular forms by:
\begin{gather}
  \rmU_{a,b}\,\sum_{n \in \ZZ} c\big(f;\,n;\,\Im(\tau)\big) e(n \tau)
:=
  \sum_{\substack{n \in \ZZ\\n \equiv b \,\pmod{a}}}
  c\big( f;\,n;\,\tfrac{\Im(\tau)}{a} \big)
  e\big( \tfrac{n \tau}{a} \big)
\tx{.}
\end{gather}

We now present a modular substitution formula for harmonic Maass forms.  We follow the notation of~\cite{radu-2013}, since our proof is similar to the proof of Lemma 5.1 of~\cite{radu-2013}, in which we replace the~$\eta$\nbd{}multiplier with $\nu_{\theta}$. 

\begin{lemma}
\label{la:qexpansion}
Let $F(\tau) = \sum_{n \in \ZZ} a(n,\Im(\tau))\, e(n \tau) \in \bbM_{k \slash 2} (\Gamma_0(4N) )$ with odd $N, k\in\ZZ_{\ge 1}$. Let $m\in \ZZ_{\ge 1}$, $t\in\{0,\ldots,m-1\}$, and~$\gamma = \begin{psmatrix} 1 & B \\ 4N C & D \end{psmatrix} \in \Gamma_0(4N)$ such that $C \isdiv m$ and $\gcd(m \slash C, 4NC) = 1$. Then
\begin{gather}
\label{eq:qexpansion}
  \big( \rmU_{m,t}\, F \bigvert_{\frac{3}{2}} \gamma \big) (\tau)
=
  \mfrac{C}{m}
  \sum_{d \isdiv \frac{m}{C}}
  e\big( \mfrac{xt}{m} \big)\,
  \epsilon_d^{-k} \delta_d\,
  d^{\frac{k}{2}}
  \sum_{n \in \ZZ}
  a\big( Cn + \ov{d}^2t, \tfrac{d^2}{m}\Im(\tau)\big) T(n,d)\,
  e\Big( \mfrac{d^2 (Cn+\ov{d}^2t)}{m}\, \tau \Big)
\tx{,}
\end{gather}
where~$\ov{d}$ is the smallest non-negative integer such that $d \ov{d} \equiv 1 \,\pmod{C}$, $x,y \in \ZZ$ satisfy
\begin{gather*}
  4 N C x + \mfrac{m}{C} y = 1
\tx{,}
\end{gather*}
and $x \equiv 0 \pmod{C}$,
\begin{gather*}
  T(n,d)
:=
  e\big( \tfrac{d^2 Dn x}{m/C} \big)
  \left(\mfrac{NC}{m \slash (Cd)} \right)
  \sum_{\substack{ 0 \le s < \frac{m}{Cd} \\ \gcd(s, \frac{m}{Cd}) = 1}}
  \left(\mfrac{s}{m \slash (Cd)} \right)
  e\Big( \mfrac{ -(x \slash C) (\ov{s} (\ov{d}^2 t + Cn) + st)}{m \slash (Cd)} \Big)
\tx{,}
\end{gather*}
where $\ov{s} s \equiv 1 \pmod{\frac{m}{Cd}}$, and
\begin{gather*}
  \delta_d
:=
  (-1)^{\frac{d-1}{2} \cdot \frac{(Nm/Cd) - 1}{2}}
  \left( \mfrac{NC}{m \slash (Cd)} \right)\,
  \left( \mfrac{d}{N} \right)\,
\tx{.}
\end{gather*}
\end{lemma}

\begin{proof}
We first expand $\rmU_{m,t}\, F$ as a sum:
\begin{gather*}
  \big( \rmU_{m,t}\, F \bigvert_{\frac{k}{2}}\, \gamma \big) (\tau)
=
  m^{\frac{k}{4} - 1} \sum_{u \pmod{m}}
  e\big( \tfrac{-ut}{m} \big)
  \Big( F \bigvert_{\frac{k}{2}}\, \begin{psmatrix} 1 & u \\ 0 & m \end{psmatrix} \gamma \Big)(\tau)
\tx{.}
\end{gather*}
Set $d = \gcd(1 + 4NCu, m)$. The fact that~$F$ is modular with respect to~$\Gamma_0(4N)$ gives:
\begin{gather*}
  F \bigvert_{\frac{k}{2}}
  \begin{psmatrix} 1 & u \\ 0 & m \end{psmatrix} \gamma
=
  F \bigvert_{\frac{k}{2}}
  \begin{psmatrix} (1+4NCu) \slash d & - Y \\ 4NC m \slash d & X \end{psmatrix}
  \begin{psmatrix} d & (B+Du)X + DmY \\ 0 & m \slash d \end{psmatrix}
=
  \epsilon_X^{-k}
  \left( \mfrac{4NCm \slash d}{X} \right)
  F \bigvert_{\frac{k}{2}}
  \begin{psmatrix} d & (B+Du)X + DmY \\ 0 & m \slash d \end{psmatrix}
\tx{,}
\end{gather*}
where $X,Y \in \ZZ$ satisfy
\begin{gather*}
  \mfrac{1+4NCu}{d}\, X + \mfrac{4NCm}{d}\, Y
=
  1
\tx{.}
\end{gather*}

Observe that $X \equiv d \,\pmod{4NC}$, so that $\epsilon_d = \epsilon_X$. Properties of the Jacobi symbol yield that:
\begin{align*}
  \left( \mfrac{4NCm \slash d}{X} \right)
&{}=
  \left( \mfrac{N \big(m \slash (Cd)\big)}{X} \right)
=
  (-1)^{\frac{d-1}{2} \cdot \frac{(Nm \slash Cd) - 1}{2}}
  \left( \mfrac{ X}{N \big(m \slash (Cd)\big)} \right)
\\
&{}=
  (-1)^{\frac{d-1}{2} \cdot \frac{(Nm \slash Cd) - 1}{2}}
  \left( \mfrac{ (1 + 4NCu) \slash d}{m \slash (Cd)} \right)\,
  \left( \mfrac{d}{N} \right)
=
  \delta_d\,
  \left( \mfrac{s}{m \slash (Cd)} \right)
\tx{,}
\end{align*}
where the last equality follows after the substitution $u=-x + sd + \frac{m}{C}r$, which is justified by~Lemma~\ref{la:residueclasses}.

We now return to $\rmU_{m,t}\, F \bigvert_{k \slash 2} \gamma$. We use again Lemma~\ref{la:residueclasses} and replace $u = -x + sd + \frac{m}{C} r$ to discover that:
\begin{gather}
\label{eq:Umt_F_rewrite}
\begin{aligned}
&
  \rmU_{m,t}\, F \bigvert_{\frac{k}{2}}\gamma
=
  m^{\frac{k}{4} - 1}
  e\big( \tfrac{xt}{m} \big)
  \sum_{d \isdiv \frac{m}{C}}
  \sum_{\substack{ 0 \le s < \frac{m}{Cd} \\ \gcd(s, \frac{m}{Cd}) = 1}}
  e\big( \tfrac{-sdt}{m} \big)
  \sum_{r=0}^{C-1}
  e\big( \tfrac{-rt}{C} \big)
  F \bigvert_{\frac{k}{2}}
  \begin{psmatrix} 1 & u \\ 0 & m \end{psmatrix} \gamma
\\
={}&
  m^{\frac{k}{4} - 1}
  e\big( \tfrac{xt}{m} \big)
  \sum_{d \isdiv \frac{m}{C}}
  \sum_{\substack{ 0 \le s < \frac{m}{Cd} \\ \gcd(s, \frac{m}{Cd}) = 1}}
  e\big( \tfrac{-sdt}{m} \big)\,
  \epsilon_d^{-k} \delta_d\,
  \left( \mfrac{s}{m \slash (Cd)} \right)
  \sum_{r=0}^{C-1}
  e\big( \tfrac{-rt}{C} \big)
  F \bigvert_{\frac{k}{2}} \begin{psmatrix} d & (B+Du)X + DmY \\ 0 & m \slash d \end{psmatrix} 
\tx{.}
\end{aligned}
\end{gather}
Next we rewrite the inner sum over~$r$ in~\eqref{eq:Umt_F_rewrite}. Observe~\eqref{eq:slash_k} and also that $F$ is translation invariant to find that: 
\begin{align}
\label{eq:inner_sum_over_r}
  \sum_{r=0}^{C-1}
   e\big( \tfrac{-rt}{C} \big)\,
   \Big(
   F \bigvert_{\frac{k}{2}}
   \begin{psmatrix} d & (B+Du)X + DmY \\ 0 & m \slash d \end{psmatrix}
   \Big)(\tau)
=
  m^{-\frac{k}{4}} d^{\frac{k}{2}}\,
  \sum_{r=0}^{C-1}
  e\big(\mfrac{-rt}{C}\big)
  F \Big( \mfrac{d\tau + (B+Du)X}{m \slash d} \Big)
\tx{.}
\end{align}

We now want to confirm that
\begin{gather}%
\label{eq:la:qexpansion:main_congruence}
  (B+Du) X
\equiv
  \mfrac{m}{C} (Drd + Byd + Dsd^2 y) + x (-\ov{s} x + Dd)
  \;\pmod{\tfrac{m}{d}}
\tx{.}
\end{gather}
Since~$\gcd(C, \frac{m}{Cd}) = 1$ and~$x \equiv 0 \,\pmod{C}$, the Chinese remainder theorem reduces this to verifying that
\begin{alignat*}{2}
  (B+Du) X
&{}\equiv
  \mfrac{m}{C} (Drd + Byd + Dsd^2 y)
  &&\;\pmod{C}
\quad\tx{and}\\
  (B+Du) X
&{}\equiv
  x (-\ov{s} x + Dd)
  &&\;\pmod{\tfrac{m}{Cd}}
\tx{.} 
\end{alignat*}
The first congruence follows from~$D \equiv 1 \,\pmod{C}$, $X \equiv d \,\pmod{C}$, $u \equiv sd + \frac{m}{C}{r} \,\pmod{C}$, and~$\frac{m}{C} y \equiv 1 \,\pmod{C}$. The second one follows from~$D = 1 + 4 N C B$, $u \equiv -x + sd \,\pmod{\frac{m}{Cd}}$, $4 NC x \equiv 1 \,\pmod{\frac{m}{Cd}}$, and~$\frac{1 + 4 NC u}{d} X \equiv 1 \,\pmod{\frac{m}{Cd}}$, which simplifies to~$4NC s X \equiv 1 \,\pmod{\frac{m}{Cd}}$.

Congruence~\eqref{eq:la:qexpansion:main_congruence} together with the facts that~$D \equiv 1 \,\pmod{C}$ and that~$F$ is translation invariant imply
\begin{gather*}
  F \Big( \mfrac{d\tau + (B+Du)X}{m \slash d} \Big)
=
  F \Big( \mfrac{ d \tau + \frac{m}{C} (Drd + Byd + Dsd^2 y) + x (-\ov{s} x + Dd)}{m \slash d} \Big)
=  F \Big( \mfrac{\tau'(s,d)+ d^2 r}{C} \Big)
\tx{,}
\end{gather*}
where
\begin{gather*}
  \tau'(s,d)
:=
  \mfrac{ d \tau + x (-\ov{s} x + Dd)}{m \slash (Cd)} +  d^2y( B + sd )
\tx{.}
\end{gather*}

Thus, for the sum over~$r$ in \eqref{eq:inner_sum_over_r} we have:
\begin{gather}
\label{eq:inner_sum_over_r_rewritten}
\begin{aligned}
  \sum_{r=0}^{C-1}
  e\big( \tfrac{-rt}{C} \big)\,
  F \bigvert_{\frac{k}{2}}
  \begin{psmatrix} d & (B+Du)X + DmY \\ 0 & m \slash d \end{psmatrix}
&{}=
  m^{-\frac{k}{4}} d^{\frac{k}{2}}\,
  \sum_{r=0}^{C-1}
  e\big( \tfrac{-r d^2 (t\ov{d}^2)}{C} \big)\,
  F \Big( \mfrac{\tau' + d^2 r}{C} \Big)
\\
&{}=
  m^{-\frac{k}{4}} d^{\frac{k}{2}} C\,
  \sum_{n \in \mathbb{Z}}
  a\big( Cn + \ov{d}^2 t, \tfrac{d^2}{m}\Im(\tau) \big)\,
  e\Big( \mfrac{\tau' (Cn + \ov{d}^2 t)}{C} \Big)
\tx{.}
\end{aligned}
\end{gather}
Consider the exponential term of the last equation. We substitute~$y \slash C=(1-4 N C x) \slash m$ and note that~$d^2 y (B + sd) Cn \slash C$ and~$st 4 NC x d \slash m$ are integers to find that
\begin{gather}
\label{eq:rewrite_exponential_term}
\begin{aligned}
&
  e\Big( \mfrac{\tau' (Cn + \ov{d}^2 t)}{C} \Big)
=
  e\Big( \mfrac{x(-\ov{s} x + D d)(C n + \ov{d}^2t)}{m \slash d} \Big)\,
  e\Big( \mfrac{d^2 y (B + sd)(C n + \ov{d}^2t)}{C} \Big)\;
  e\Big( \mfrac{d^2 (Cn + \ov{d}^2 t)}{m} \tau \Big)
\\
={}&
  e\Big( \mfrac{tyB}{C} + \mfrac{d^2 x D \ov{d}^2 t}{m} \Big)\,
  e\Big( \mfrac{d^2 Dn x}{m \slash C} \Big)\,
  e\Big( \mfrac{ - x^2 \ov{s} (Cn + \ov{d}^2 t)}{m \slash d} \Big)\,
  e\Big( \mfrac{s d t}{m} \Big)\;
  e\Big( \mfrac{d^2 (C n + \ov{d}^2 t)}{m} \tau \Big)
\tx{.}
\end{aligned}
\end{gather}

We now go back to $\rmU_{m,t} F |_{\frac{k}{2}} \gamma (\tau)$ in~\eqref{eq:Umt_F_rewrite}, and use~\eqref{eq:inner_sum_over_r_rewritten} and~\eqref{eq:rewrite_exponential_term} to we rewrite the two inner sums: 
\begin{align*}
&
  \sum_{\substack{ 0 \le s < \frac{m}{Cd} \\ \gcd(s, \frac{m}{Cd}) = 1}}
  e\big(\tfrac{-sdt}{m}\big)
  \sum_{r=0}^{C-1}
  e\big(\tfrac{-rt}{C}\big)
  F \bigvert_{\frac{k}{2}}
  \begin{psmatrix} 1 & -x + sd + r m \slash C \\ 0 & m \end{psmatrix} \gamma
\\
={}&
  m^{-\frac{k}{4}} d^{\frac{k}{2}} C\,
  \sum_{\substack{ 0 \le s < \frac{m}{Cd} \\ \gcd(s, \frac{m}{Cd}) = 1}}
  e\big(\tfrac{-sdt}{m}\big)\,
  \epsilon_d^{-k} \delta_d\,
  \left( \mfrac{s}{m \slash (Cd)} \right)\,
  \sum_{n \in \mathbb{Z}}
  a\big( Cn + \ov{d}^2 t, \tfrac{d^2}{m}\Im(\tau) \big)
  e\Big( \mfrac{\tau' (Cn + \ov{d}^2 t)}{C} \Big)
\\
={}&
  m^{-\frac{k}{4}} d^{\frac{k}{2}} C\,
  \epsilon_d^{-k} \delta_d\,
  e\big( \tfrac{tyB}{C} + \tfrac{d^2 x D \ov{d}^2 t}{m} \big)
  \sum_{n \in \ZZ}
  a\big( Cn + \ov{d}^2 t, \tfrac{d^2}{m}\Im(\tau)\big)
  T'(n,d)\,
  e\Big( \mfrac{d^2 (Cn+\ov{d}^2t)}{m}\, \tau \Big)
\tx{,}
\end{align*}
where
\begin{gather*}
  T'(n,d)
=
  e\big( \tfrac{d^2 Dn x}{m \slash C} \big)\,
  \sum_{\substack{ 0 \le s < \frac{m}{Cd} \\ \gcd(s, \frac{m}{Cd}) = 1}}
  \left( \mfrac{s}{m \slash Cd} \right)
  e\big( \tfrac{ - x^2 \ov{s} (Cn + \ov{d}^2 t)}{m \slash d} \big)
  e\big( \tfrac{- s t 4 N C x}{m \slash d} \big)
\tx{.}
\end{gather*}

Now, recall that~$\gcd(4 N C, m \slash C) = 1$. Hence $s\mapsto 4NCs$ is a bijection modulo $m\slash (Cd)$. That substitution together with replacing $4NCx$ by $1-(m\slash C)y$ then yields that~$T'(n,d) = T(n,d)$, which completes the proof.
\end{proof}

The following lemma was used in the proof of Lemma~\ref{la:qexpansion}, and it is a minor modification of an argument from~\cite{radu-2013}.
\begin{lemma}%
\label{la:residueclasses}
Let $C,m,d \in \ZZ_{\ge 1}$ and $\alpha, x, y \in \ZZ$. Suppose that $C \isdiv m$, $\gcd(\alpha C, \frac{m}{C}) = 1$, $d \isdiv \frac{m}{C}$ and $\alpha Cx + \frac{m}{C}y= 1$. Then 
\begin{gather*}
  \big\{ u \,\pmod{m} \,:\, \gcd(m, 1 + \alpha Cu) = d \big\}
=
  \big\{ -x + sd + \tfrac{m}{C}r \,\pmod{m} \,:\, \gcd(s, \tfrac{m}{Cd}) = 1 \big\}
\tx{.}
\end{gather*}
\end{lemma}

\begin{proof}
Given an element~$-x + sd + \frac{m}{C} r$ of the right hand side, we observe that
\begin{align*}
  \gcd\big(m, 1 + \alpha C(-x + sd + \tfrac{m}{C}r)\big)
&{}=
  \gcd\big(m, 1 - \alpha Cx + \alpha Csd\big)
=
  \gcd\big(m, \tfrac{m}{C}y + \alpha Csd\big)
\\
&{}=
  d\,
  \gcd(\tfrac{m}{d}, \tfrac{m}{Cd}y + \alpha Cs)
=
 d
\tx{,}
\end{align*}
where we used the property~$\gcd(s, \frac{m}{Cd}) = 1$ in the last equality. This shows that the right hand side is contained in the left hand side.

Now consider an element~$u$ of the left hand side. Let~$s = \frac{(1 + \alpha Cu)x}{d}$, which is an integer since~$\gcd(m, 1+\alpha C u) = d$. Then $-x + sd = \alpha xCu \equiv u \,\pmod{\frac{m}{C}}$, i.e., $u = -x + sd + \frac{m}{C}r$ for some $r \in \ZZ$. This shows that the left hand side is contained in the right hand side.

We observe that each class on the right hand side represents only one class on the left hand side with the restrictions $0 \le s < \frac{m}{Cd}, 0 \le r < C$. First, from
\begin{gather*}
  -x + s_1d + \tfrac{m}{C}r_1
\equiv
  -x + s_2d + \tfrac{m}{C}r_2
  \;\pmod{m}
\end{gather*}
we conclude that~$s_1 - s_2 \equiv \frac{m}{Cd}(r_2 - r_1) \,\pmod{\tfrac{m}{d}}$.
This implies that $s_1 - s_2 \equiv 0 \pmod{\frac{m}{Cd}}$, and since $0 \le s_1, s_2 < \frac{m}{Cd}$ we must have $s_1 = s_2$. 
Hence $\frac{m}{Cd} (r_2 - r_1) \equiv 0\pmod{\frac{m}{Cd}C}$, which implies (since $\gcd(\frac{m}{Cd}, C ) = 1$) that~$r_1 - r_2 \equiv 0 \pmod{C}$, and hence $r_1 = r_2$. 
\end{proof}

\section{Proof of Theorem~\ref{mainthm:ramanujan_type_and_discriminant_conditions}}
\label{sec:Proof_of_Theorem_B}

We start this section with a proposition that is a weaker version of Theorem~\ref{mainthm:ramanujan_type_and_discriminant_conditions}, where the analog of the set~$S_0$ in Theorem~\ref{mainthm:ramanujan_type_and_discriminant_conditions} is subdivided into two sets~$S_{0+}$ and~$S_{0-}$. 

\begin{proposition}
\label{prop:ramanujan_type_and_discriminant_conditions_S0pm}
Let~$\ell \ge 3$ be prime.  Assume that there is a maximal regular Ramanujan-type congruence for Hurwitz class numbers modulo~$\ell$ on~$a \ZZ + b$ for cube-free, odd~$a$. Define a disjoint union of finite sets of odd primes~$S = S_+ \cup S_- \cup S_{0+} \cup S_{0-}$ via the formulae
\begin{alignat*}{2}
  S_+ &{}= \big\{ p \,:\, p \isexactdiv a, \big(\tfrac{-b}{p}\big) = +1 \big\}
\tx{,}\quad
&
  S_- &{}= \big\{ p \,:\, p \isexactdiv a, \big(\tfrac{-b}{p}\big) = -1 \big\}
\tx{,}
\\
  S_{0+} &{}= \big\{ p \,:\, p^2 \isexactdiv a, \big(\tfrac{-b}{p}\big) = 0, \big(\tfrac{-b \slash p}{p}\big) = +1 \big\}
\tx{,}\quad
&
  S_{0-} &{}= \big\{ p \,:\, p^2 \isexactdiv a, \big(\tfrac{-b}{p}\big) = 0, \big(\tfrac{-b \slash p}{p}\big) = -1 \big\}
\tx{,}
\end{alignat*}
where~$p$ always denotes an odd prime. Then for any imaginary quadratic field~$L$ that
\begin{enumerateroman}
\item is split at each prime of\/~$S_+$,
\item is inert at each prime of\/~$S_-$,
\item is ramified at each prime of\/~$S_{0+} \cup S_{0-}$,
\item whose discriminant~$-D_L$ of~$L$ satisfies~$\big( \frac{-D_L \slash p}{p} \big) = \pm 1$ for~$p \in S_{0\pm}$,
\end{enumerateroman}
the class number~$h_L$ is divisible by~$\ell$.

Conversely, if the conclusion in the previous paragraph holds for some disjoint union of finite sets of odd primes~$S = S_+ \cup S_- \cup S_{0+} \cup S_{0-}$, then there is a maximal regular Ramanujan-type congruence for Hurwitz class numbers modulo~$\ell$ on all~$a \ZZ + b$ with
\begin{gather*}
  a = \prod_{p \in S_+ \cup S_-} p\,
      \prod_{p \in S_{0+} \cup S_{0-}} p^2
\tx{,}\\
  b \in \ZZ \tx{\ such that\ }
  \big( \tfrac{-b}{p} \big) = \pm 1 \tx{\ if\ } p \in S_\pm
  \tx{,\ }
  p \isexactdiv b \tx{\ if\ } p \in S_{0+} \cup S_{0-}
  \tx{\ and\ }
  \big( \tfrac{-b \slash p}{p} \big) = \pm 1 \tx{\ if\ } p \in S_{0\pm}
\tx{.}
\end{gather*}
\end{proposition}

\begin{proof}
Assume that there is a maximal regular Ramanujan-type congruence for Hurwitz class numbers modulo~$\ell$ on~$a \ZZ + b$, which by Corollary~\ref{cor:aktin_type_ramanujan_type} is then also an Atkin-type congruence for Hurwitz class numbers modulo~$\ell$ around~$b \pmod{a}$. Note that~$S$ contains all prime divisors of~$a$. Now, if~$L$ is an imaginary quadratic field with discriminant~$-D_L$ as in the statement, then the assumptions on~$L$ and~$-D_L$ imply that~$D_L \in b (\ZZ \slash a \ZZ)^{\times\,2}$. Thus, $H(D_L) \equiv 0 \,\pmod{\ell}$ and hence $h(-D_L) = H(D_L) \omega(-D_L) \slash 2 \equiv 0 \,\pmod{\ell}$, where~$\omega(-D_L)$ is as in~\eqref{eq:hurwitz_class_number}.

Next we consider the converse. If~$D f^2 \in b (\ZZ \slash a \ZZ)^{\times\,2}$ with a negative fundamental discriminant~$-D$, then the conditions on~$-b$ and the fact that~$\ord_p(a) > \ord_p(b)$ for all~$p \in S$ imply that
\begin{gather*}
  \big(\tfrac{-D}{p}\big) = \pm 1
  \tx{\ if\ } p \in S_\pm
\quad\tx{and}\quad
  p \isdiv -D
  \tx{\ as well as\ }
  \big(\tfrac{-D \slash p}{p}\big) = \pm 1
  \tx{\ if\ } p \in S_{0\pm}
\tx{.}
\end{gather*}
By assumption, $h(-D)\equiv0\,\pmod{\ell}$. Therefore, $-D \not\in \{-3,-4\}$, $H(D)=h(-D)\equiv0\, \pmod{\ell}$, and~\eqref{eq:hurwitz_class_number} shows that~$H(D f^2) \equiv 0 \,\pmod{\ell}$.
\end{proof}

The remaining task in this section is to establish a relation between the sets~$S_{0+}$ and~$S_{0-}$ in~Proposition~\ref{prop:ramanujan_type_and_discriminant_conditions_S0pm}.  The next theorem shows that congruences associated with primes in~$S_{0\pm}$ are ``coupled'' to each other in the sense that they apply one another.

\begin{theorem}
\label{thm:congruences_coupled_square_classes}
Let~$\ell \ge 3$ be prime.  Assume that there is a Ramanujan-type congruence for Hurwitz class numbers modulo~$\ell$ on~$a p^2 \ZZ + b$ for a positive integer~$a$, some prime~$p$, and an integer~$b$. Furthermore, assume that~$p \nisdiv a$, $p \ne \ell$, and~$p \isexactdiv b$. Then the Hurwitz class numbers have a Ramanujan-type congruence modulo~$\ell$ on~$a p^2 \ZZ + b'$ for all integers~$b'$ such that~$p \isexactdiv b'$ and~$b \equiv b' \,\pmod{a}$.
\end{theorem}

\begin{proof}%
Recall Zagier's~\cite{zagier-1975} Eisenstein series:
\begin{gather}
\label{eq:zagier-eisenstein-series}
  E_{\frac{3}{2}}(\tau)
\;:=\;
  \sum_{D = 0}^\infty H(D) e(D \tau)
\,+\,
  \mfrac{1}{16 \pi} \theta^\ast(\tau)
\;\in\;
  \bbM_{\frac{3}{2}}(\Gamma_0(4))
\tx{,}
\end{gather}
where $\theta^\ast$ is the non-holomorphic Eichler integral of a classical Jacobi theta constant (see for example~\cite{beckwith-raum-richter-2020,beckwith-raum-richter-2022}).  The assumption that~$p$ exactly divides~$b$ guarantees that $-b$ is not a square modulo~$a p^2$ and hence $\rmU_{a p^2, b} E_{\frac{3}{2}}$ is a holomorphic modular form of level $4a p^2$.  The assumption that~$H(a p^2 n + b) \equiv 0 \pmod{\ell}$ then gives~$\rmU_{a p^2, b} E_{\frac{3}{2}} \equiv 0 \pmod{\ell}$. The main theorem of~\cite{beckwith-raum-richter-2020} asserts that~$\ell \isdiv a$. Thus, Proposition~2.2 of~\cite{ahlgren-beckwith-raum-2020-preprint} (a version of the $q$-expansion principle) then implies that~$\rmU_{a p^2, b}\, E_{\frac{3}{2}} |_{\frac{3}{2}} \gamma \equiv 0 \pmod{\ell}$ where $\gamma = \begin{psmatrix} 1 & \ast \\ 4a & D \end{psmatrix} \in \Gamma_0(4)$, where~$D$ is a suitable integer.

Apply Lemma~\ref{la:qexpansion}, with $F = E_{\frac{3}{2}}$ (so that $N=1$), $m = a p^2$, $t = b$, and $C = a$. Then the formula for $\rmU_{a p^2, b} E_{\frac{3}{2}} |_{\frac{3}{2}} \gamma$ in~\eqref{eq:qexpansion} is a double sum over $d | p^2$ and over $n\in\ZZ$. The terms in that double sum for which the exponent of $q$ is in $\frac{1}{a p} \mathbb{Z} \backslash \frac{1}{a} \mathbb{Z}$ all arise from the terms with $d=1$, and these terms must sum to~$0$ modulo~$\ell$:
\begin{gather*}
  \sum_{p \isexactdiv (a n + b)}
  H(a n + b) T(n,1)\,
  e\big( \tfrac{a n + b}{a p^2} \tau \big)
\equiv
  0
  \;\pmod{\ell}
\tx{.}
\end{gather*}

If~$p\isexactdiv  (an + b)$ and~$\ov{4a}$ and~$\ov{s}$ are integers such that~$4a \ov{4a} \equiv 1 \pmod{p^2}$ and~$s \ov{s} \equiv 1 \pmod{p^2}$, then 
\begin{align*}
  T(n,1)
&=
  e\big( \tfrac{Dn \ov{4a}}{p^2} \big)
  \left(\frac{a}{p^2} \right)
  \sum_{\substack{0 \le s < p^2 \\ p \nmid{s}}}
  \left( \frac{s}{p^2} \right)
  e\Big(  \frac{-\frac{\ov{4a}}{a} (\ov{s} (b + a n) + sb)}{p^2} \Big)
\\
&=
  e\big( \tfrac{Dn \ov{4a}}{p^2} \big)
  K\big( -\tfrac{\ov{4a}b}{a}, \tfrac{-\ov{4a}(b + an)}{a}, p^2 \big)
=
  e\big( \tfrac{Dn \ov{4a}}{p^2} \big)
  p
  K\big( -\tfrac{\ov{4a}b}{a p}, \tfrac{-\ov{4a}(b + an)}{a p}, p \big)
\tx{,}
\end{align*}
where
\begin{gather*}
K(a',b',c')
=
  \sum_{\substack{0 \le s < c' \\ \gcd(s,c')=1}}
  e\big( \tfrac{a's + b'\ov{s}}{c'} \big)
\tx{,}
\end{gather*}
with~$\ov{s}$ a multiplicative inverse of~$s$ modulo~$c'$, denotes the usual Kloosterman sum for integers $a'$, $b'$, and $c'$.  Since $p\not=\ell$, $p \isdiv b$, and $p \isdiv  (an + b)$,  Lemma~5.2 of~\cite{ahlgren-beckwith-raum-2020-preprint} implies that $T(n,1) \not\equiv 0 \pmod{\ell}$. Thus, if~$p\isexactdiv  (an + b)$, then 
\begin{gather*}
  H(an + b) \equiv 0 \;\pmod{\ell}
\tx{,}
\end{gather*}
which yields the claim.
\end{proof}

We are now in a position to prove Theorem~\ref{mainthm:ramanujan_type_and_discriminant_conditions}.
\begin{proof}%
[Proof of Theorem~\ref{mainthm:ramanujan_type_and_discriminant_conditions}]
We first show the converse direction. Assume that for all imaginary quadratic fields~$L$ with conditions (i), (ii), and (iii) it follows that~$h_L = h(-D_L)$ is divisible by~$\ell$. Set
\begin{gather*}
  S_{0+} = \big\{ p \in S_0 \,:\, \big(\tfrac{-D_L \slash p}{p}\big) = +1 \big\}
\quad\tx{and}\quad
  S_{0-} = \big\{ p \in S_0 \,:\, \big(\tfrac{-D_L \slash p}{p}\big) = -1 \big\}
\tx{.}
\end{gather*}
Observe that~$S_0 = S_{0+} \cup S_{0-}$.  Moreover, conditions (i), (ii), and (iii) of Proposition~\ref{prop:ramanujan_type_and_discriminant_conditions_S0pm} coincide with (i), (ii), and (iii) of Theorem~\ref{mainthm:ramanujan_type_and_discriminant_conditions}, while condition~(iv) of Proposition~\ref{prop:ramanujan_type_and_discriminant_conditions_S0pm} holds trivially.  In particular, for all imaginary quadratic fields~$L$ with conditions (i)-(iv) in Proposition~\ref{prop:ramanujan_type_and_discriminant_conditions_S0pm} it follows that~$h_L$ is divisible by~$\ell$. Hence Proposition~\ref{prop:ramanujan_type_and_discriminant_conditions_S0pm} shows that the Hurwitz class numbers have a maximal regular Ramanujan-type congruence modulo~$\ell$ on~$a \ZZ + b$ as stated in Theorem~\ref{mainthm:ramanujan_type_and_discriminant_conditions}. This proves the converse direction of the theorem.

For the other direction, assume that the Hurwitz class numbers have a maximal regular Ra\-ma\-nu\-jan-type congruence modulo~$\ell$ on~$a \ZZ + b$ and let~$S = S_+ \cup S_- \cup S_0$ be as in the first part of the theorem. Let~$S_{0+}$ and~$S_{0-}$ be as in Proposition~\ref{prop:ramanujan_type_and_discriminant_conditions_S0pm} and observe again that~$S_0 = S_{0+} \cup S_{0-}$, since~$a$ is cube-free and hence~$(\tfrac{-b \slash p}{p}) \ne 0$ for all~$p \in S_0$.

For a given imaginary quadratic field~$L$ set
\begin{gather*}
  S'_{0+} = \big\{ p \in S_0 \,:\, \big(\tfrac{-D_L \slash p}{p}\big) = +1 \big\}
\quad\tx{and}\quad
  S'_{0-} = \big\{ p \in S_0 \,:\, \big(\tfrac{-D_L \slash p}{p}\big) = -1 \big\}
\tx{.}
\end{gather*}
Then~$S_0 = S_{0+} \cup S_{0-} = S'_{0+} \cup S'_{0-}$. We need to demonstrate that~$\ell \isdiv h_L$, which will follow from Proposition~\ref{prop:ramanujan_type_and_discriminant_conditions_S0pm} if we can show that there is a Ramanujan-type congruence for Hurwitz class numbers modulo~$\ell$ on~$a \ZZ + b'$ with~$(\frac{b'}{p}) = \pm 1$ for~$p \in S_\pm$ and~$p \isdiv b'$, $(\frac{b' \slash p}{p}) = \pm 1$ for all~$p \in S'_{0\pm}$. The remainder of the proof is dedicated to showing the existence of such a Ramanujan-type congruence.

We induct on the cardinality of the symmetric difference~$S_{0+} \symdiff S'_{0+}$. If it is empty, then we have~$S_{0-} = S'_{0-}$ and the Ramanujan-type congruence on~$a \ZZ + b'$ follows from the one on~$a \ZZ + b$ by Theorem~\ref{thm:rc_on_square_classes_ell3}. If~$S_{0+} \symdiff S'_{0+} \ne \emptyset$, then there is a prime~$q \in S_{0+} \setminus S'_{0+}$ or~$q \in S_{0-} \setminus S'_{0-}$. Without loss of generality assume the former, and set~$S''_{0+} = S'_{0+} \cup \{q\}$ and~$S''_{0-} = S'_{0-} \setminus \{q\}$. Since~$\#( S_{0+} \symdiff S'_{0+} ) = \#( S_{0+} \symdiff S''_{0+} ) + 1$, we can apply the induction hypothesis to obtain a Ramanujan-type congruence on~$a \ZZ + b''$ for any integer~$b''$ with $(\frac{-b'' \slash p}{p}) = \pm 1$ for~$p \in S''_{0\pm}$. In particular, we may choose some~$b'' \equiv b \,\pmod{a_q^\#}$. For the induction step it suffices to apply Theorem~\ref{thm:congruences_coupled_square_classes} with~$b$ and~$b'$ replaced by~$b''$ and~$b$, respectively.  This completes the proof.
\end{proof}

\section{Proofs of Theorems~\ref{mainthm:indivisibility} and~\ref{mainthm:classification_nonhol_nondeg}}
\label{sec:Proof_of_Theorems_A_C}

\begin{proof}%
[Proof of Theorem~\ref{mainthm:classification_nonhol_nondeg}]
We use the Hurwitz class number formula~\eqref{eq:hurwitz_class_number}, similarly to an argument of our previous work~\cite{beckwith-raum-richter-2020}, to show that there is a Ramanujan-type congruence on $\ell^{m+1}\ZZ + u \ell^m$ of the given form: Since~$m > 0$ is even and~$\ell$ is odd, any~$D f^2 \in \ell^m (\ZZ \ell + u)$ for a negative fundamental discriminant~$-D$ satisfies~$(\frac{-D}{\ell}) = (\frac{-u}{\ell}) = +1$. In particular, if~$\ell = 3$, for which~$H(D)$ is not evidently~$\ell$\nbd{}integral, then~$3 \nisdiv -D$, that is~$-D \ne -3$, and hence~$H(D)$ is $3$-integral.  Now the contribution of~$p = \ell$ in the product over~$p \isdiv f$ in~\eqref{eq:hurwitz_class_number} is given by
\begin{gather}
\label{eq:mainthm:classification_nonhol_nondeg_hecke_factor}
  1 \,+\,
  \sum_{\substack{d \isdiv f\\ d \tx{\ an\ }\ell\tx{-power}\\ d \ne 1}} \mspace{-20mu}
  \big( d - \tfrac{d}{\ell} \big(\tfrac{-D}{\ell}\big) \big)
\equiv
  1 - \big(\tfrac{-D}{\ell}\big)
  \;\pmod{\ell}
\tx{,}
\end{gather}
which vanishes modulo~$\ell$.  The resulting Ramanujan-type congruence is non-holomorphic, since~$-\ell^m u$ is a square modulo~$\ell^{m+1}$ by the assumptions on~$m$ and~$u$.

To demonstrate that every maximal regular, non-holomorphic Ramanujan-type congruence is of the given form, recall from Section~\ref{sec:atkin_type_congruence} that Ramanujan-type and Atkin-type congruences for Hurwitz class numbers imply each other. Thus, we can assume that we have a maximal regular, non-holomorphic Atkin-type congruence of the Hurwitz class numbers modulo~$\ell$ around~$b \pmod{a}$. Then~$\ell \isdiv a, b$ by Theorem~\ref{thm:ell_div_ab_atkin_type}.  The congruence is regular by assumption and hence~$\ell^2 \isdiv a$.  Since the congruence is non-holomorphic, it is impossible that~$\ell \isexactdiv b$, and thus~$\ell^2 \isdiv b$. Again, since the congruence is regular, we find that~$\ell^3 \isdiv a$, which proves one of the claims of the theorem. The assumption that the Atkin-type congruence is maximal in combination with Proposition~\ref{prop:maximal_atkin_type_ordp} yield that~$a \isdiv \ell b$. Again, since the congruence is regular, the stronger statement~$a \isexactdiv \ell b$ holds, which establishes another claim of the theorem.

It remains to show that~$a$ is an odd power of~$\ell$ and that~$-\ell b \slash a$ is a square modulo~$\ell$.  To this end, we inspect the Hurwitz class number formula~\eqref{eq:hurwitz_class_number} for arbitrary elements~$D f^2 \in b (\ZZ \slash a \ZZ)^{\times\,2}$, where~$-D$ is a negative fundamental discriminant. Note that~$\ell \isdiv f$, since~$\ell^2 \isdiv b$ and~$\ell$ is odd.  Furthermore, note that~$f_\ell$ coincides for all~$D f^2 \in b (\ZZ \slash a \ZZ)^{\times\,2}$.  Assume by contraposition that the contribution of~$p = \ell$ to the product over~$p \isdiv f$ in~\eqref{eq:hurwitz_class_number} does not vanish, i.e, that the congruence~\eqref{eq:mainthm:classification_nonhol_nondeg_hecke_factor} does not vanish.  Then~\eqref{eq:hurwitz_class_number} implies an Atkin-type congruence modulo~$\ell$ around~$(b \slash f_\ell^2) \pmod{a}$. We have~$\ord_\ell(a \slash f_\ell^2) = \ord_\ell(b \slash f_\ell)^2 + 1$ and thus by Proposition~\ref{prop:maximal_atkin_type_ordp} the Atkin-type congruence around~$(b \slash f_\ell^2) \pmod{a}$ is not maximal and implies a maximal regular one around~$(b \slash f_\ell^2) \pmod{a \slash f_\ell^2}$. The matter that~$\ell^2 \nisdiv (b \slash f_\ell^2)$ is a contradiction to what we have already established. Thus, the congruence~\eqref{eq:mainthm:classification_nonhol_nondeg_hecke_factor} vanishes modulo~$\ell$.

As before in the proof of the second part of the theorem one finds that~$\ell \nisdiv -D$ and~$-D$ is a square modulo~$\ell$, and if~$\ell = 3$, then~$H(D)$ is~$3$-integral.  The vanishing of~\eqref{eq:mainthm:classification_nonhol_nondeg_hecke_factor} modulo~$\ell$ implies an Atkin-type congruence around~$b \pmod{a_\ell}$.  Since the one around~$b \pmod{a}$ is maximal regular by assumption, we have~$a = a_\ell$. Hence~$a$ is a power of~$\ell$ (with necessarily odd exponent) as claimed. Finally, the findings that~$\ell \nisdiv -D$ and that $a$ is an odd power of~$\ell$ yield that~$-D$ lies in the same square class modulo~$\ell$ as~$-\ell b \slash a = -b \slash f_\ell^2$, which shows that~$(\frac{-\ell b \slash a}{\ell}) = +1$.

\end{proof}

\begin{proof}%
[Proof of Theorem~\ref{mainthm:indivisibility}]
The case~$\ell = 2$ follows from Gauss' genus theory. Let~$\ell \ge 3$ and suppose by contradiction that for all imaginary quadratic fields~$L$ that are split at~$S_+$, $h_L$ is divisible by $\ell$. Theorem~\ref{mainthm:ramanujan_type_and_discriminant_conditions} then asserts that the Hurwitz class numbers have a maximal regular Ramanujan-type congruence modulo~$\ell$ on an arithmetic progression~$a \ZZ + b$ with~$a = \prod_{p \in S_+} p$. This contradicts Theorem~\ref{mainthm:classification_nonhol_nondeg}.
\end{proof}


\vspace{1.5\baselineskip}
\ifbool{nobiblatex}{%
  \sloppy
  \bibliographystyle{alpha}%
  \bibliography{bibliography.bib}%
}{%
  \sloppy
  \Needspace*{4em}
  \printbibliography[heading=none]
}


\Needspace*{3\baselineskip}
\noindent
\rule{\textwidth}{0.15em}

{\noindent\small
Olivia Beckwith\\
Mathematics Department,
Tulane University,
New Orleans, LA 70118, USA\\
E-mail: \url{obeckwith@tulane.edu}\\
Homepage: \url{https://www.olivia-beckwith.com/}
}\vspace{.5\baselineskip}

{\noindent\small
Martin Raum\\
Chalmers tekniska högskola och G\"oteborgs Universitet,
Institutionen för Matematiska vetenskaper,
SE-412 96 Göteborg, Sweden\\
E-mail: \url{martin@raum-brothers.eu}\\
Homepage: \url{http://raum-brothers.eu/martin}
}\vspace{.5\baselineskip}

{\noindent\small
Olav K. Richter\\
Department of Mathematics,
University of North Texas,
Denton, TX 76203, USA\\
E-mail: \url{richter@unt.edu}\\
Homepage: \url{http://www.math.unt.edu/~richter/}
}


\ifdraft{%
\listoftodos%
}

\end{document}
